\documentclass[11pt,a4paper]{article}
\usepackage[utf8]{inputenc}
\usepackage[german,english]{babel}
\usepackage{amsmath,amssymb,amsthm}
\usepackage{geometry}
\newtheorem{theorem}{Theorem}[section]
\newtheorem{lemma}[theorem]{Lemma}
\newtheorem{corollary}[theorem]{Corollary}

\theoremstyle{definition}
\newtheorem{remark}[theorem]{Remark}

\numberwithin{equation}{section}

\title{\textbf{Uniqueness results for the exterior Dirichlet and Neumann problem in two dimensions for the Euler-Lagrange equations of variational integrals with linear growth}}
\author{\textbf{Michael Bildhauer} \quad \textbf{Martin Fuchs}}
\date{}

\begin{document}
\maketitle
\begin{abstract}
We consider solutions $u, v:\Omega\rightarrow\mathbb{R}^{N}$ of the Euler-Lagrange equations associated to the
variational integral $\int_{\Omega}f(\nabla u)dx$, where $\Omega$ is an exterior domain in $\mathbb{R}^{2}$ e.g., the complement
of a disk or a half plane, and the density $f$ is of linear growth. Among other things we show
that $u|_{\partial\Omega}=v|_{\partial\Omega}$ together with the boundedness of $u-v$ yields that $u=v$ on the entire
domain.
\end{abstract}

\noindent
\textbf{Mathematics Subject Classification.} 35J25, 35J65, 49Q05, 53A10. \\
\textbf{Keywords.} exterior Dirichlet problem, exterior Neumann problem, uniqueness of solutions,
variational integrals with linear growth.

\section{Introduction and statement of the results}

In our note we investigate the exterior Dirichlet problem and in addition the exterior Neumann
problem on planar domains associated with the Euler-Lagrange equations corresponding to vari-
ational integrals $\int_{\Omega}f(\nabla u)dx$ with density $f$ being of linear growth providing some uniqueness
results. A prominent example is given by the Dirichlet problem
\begin{equation}
\begin{cases}
\text{div}\left(\frac{\nabla u}{\sqrt{1+|\nabla u|^{2}}}\right)=0 & \text{on } \Omega, \\
u=\varphi & \text{on } \partial\Omega
\end{cases}
\end{equation}
for functions $u:\Omega\rightarrow\mathbb{R}$, where $\Omega\subset\mathbb{R}^{2}$ is a domain in $\mathbb{R}^{2}$ such that $\mathbb{R}^{2}-\Omega$ is bounded and
$\varphi$ denotes a given function. Note that the first line of (1.1) is the Euler-Lagrange equation for the
area integral $\int_{\Omega}\sqrt{1+|\nabla u|^{2}}dx$, but even in the case of problem (1.1) uniqueness of solutions can
not be expected: for instance, if we let $\varphi=0$ and take $\Omega:=[|x|>1]:=\{x\in\mathbb{R}^{2}:|x|>1\}$, then
problem (1.1) admits the three different solutions $u=0$, $k_{\pm}(x):=\pm \ln(|x|+\sqrt{|x|^{2}-1})$, $|x|\ge1$,
which means that we have the trivial one together with the upper and lower catenoid. Therefore
one can only hope for uniqueness under additional restrictions to be imposed on the solutions as
it is stated in Theorem 11.3 of Osserman's monograph [1]: if $u, v$ denote bounded solutions of
(1.1), then it holds $u=v$. This means that we have uniqueness for problem (1.1) in the space
$L^{\infty}(\Omega)$, however, as Theorem 11.4 in [1] shows, the existence of solutions to (1.1) in the space
$L^{\infty}(\Omega)$ can in general not be guaranteed. The delicate question of existence (and uniqueness)
for the exterior Dirichlet problem (1.1) has been investigated by many authors leading to the
result that, roughly speaking, a smallness condition for the oscillation of the boundary data
$\varphi$ is sufficient for the existence of a solution $u$ to (1.1), whereas prescribed asymptotic growth of
$u$ at infinity or the existence of a constant asymptotic normal vector for the graph of $u$ provides
uniqueness. We refer the reader to the contributions of for instance Kuwert [2], Kutev and Tomi
[3], do Espírito-Santo and Ripoll [4], Ripoll and Tomi [5] and Aiolfi, Assmann and Ripoll [6],
where one will find further references and also a careful analysis of the historical background.

We now leave the minimal surface case and want to improve the uniqueness result from
Theorem 11.3 in [1] even in a more general framework. Our assumptions are as follows: let
$G\subset\mathbb{R}^{2}$ denote a bounded Lipschitz domain and consider the exterior region
\begin{equation}
\Omega:=\mathbb{R}^{2}-\overline{G}.
\end{equation}
For $N\ge1$ we denote by $\mathbb{R}^{2N}$ the space of $(2\times N)$-matrices. Suppose that $f:\mathbb{R}^{2N}\rightarrow\mathbb{R}$ is of
class $C^{2}$ and in addition of linear growth in the sense that
\begin{equation}
a|p|-b\le f(p)\le A|p|+B, \quad p\in\mathbb{R}^{2N},
\end{equation}
holds with positive constants $a, b, A, B$. We require further that
\begin{equation}
|D^{2}f|\in L^{\infty}(\mathbb{R}^{2N})
\end{equation}
and that
\begin{equation}
D^{2}f(p)(q,q)>0, \quad p,q\in\mathbb{R}^{2N}-\{0\}
\end{equation}
is satisfied. Due to (1.5) the density $f$ is a strictly convex function, which together with (1.3)
yields $|\nabla f|\in L^{\infty}(\mathbb{R}^{2N})$ (compare, e.g., Lemma 2.2, p.156, in [7]). We then look at weak solutions
$u\in C^{1}(\Omega,\mathbb{R}^{N})$ of the equation
\begin{equation}
\text{div}[\nabla f(\nabla u)]=0 \quad \text{on } \Omega,
\end{equation}
which by definition means that we require
\begin{equation}
\int_{\Omega}\nabla f(\nabla u):\nabla\varphi~dx=0
\end{equation}
for all $\varphi\in C_{0}^{1}(\Omega,\mathbb{R}^{N})$, ``:'' denoting the scalar product of matrices.

Our first result concerns the exterior Dirichlet problem for equation (1.6) giving an extension
of Theorem 11.3 from [1]:

\begin{theorem}
Assume that we have (1.2)--(1.5) and suppose that $u, v\in C^{1}(\Omega,\mathbb{R}^{N})\cap C^{0}(\overline{\Omega},\mathbb{R}^{N})$
are weak solutions of (1.6) in the sense of (1.7). Assume further that
\begin{equation}
u=v \quad \text{on } \partial\Omega.
\end{equation}
Then it follows $u=v$ on $\Omega$, if
\begin{equation}
u-v\in L^{\infty}(\Omega,\mathbb{R}^{N})
\end{equation}
or
\begin{equation}
\limsup_{R\rightarrow\infty}\int_{[R<|x|<2R]}|u-v|dx<\infty
\end{equation}
holds.
\end{theorem}

We note that in the case of bounded solutions $u, v$ condition (1.9) clearly holds, hence we
have uniqueness in the class of bounded solutions of the exterior Dirichlet problem associated
to equation (1.6).

From the proof of Theorem 1.1 it follows:
\begin{corollary}[Uniqueness of entire solutions]
Let (1.3)--(1.5) hold for the function $f$
and assume that $u, v\in C^{1}(\mathbb{R}^{2},\mathbb{R}^{N})$ are weak solutions of equation (1.6) on the whole plane, i.e.
(1.7) is satisfied for all $\varphi\in C_{0}^{1}(\mathbb{R}^{2},\mathbb{R}^{N})$. Then, in the case that $u-v$ is bounded, $u=v$ is true.
\end{corollary}

Next we apply Theorem 1.1 to a situation which is closer to the minimal surface setting
introduced in problem (1.1): let $N=1$ and consider a density $f$ with radial structure, i.e. we
let $f(p):=g(|p|)$, $p\in\mathbb{R}^{2}$, with $g:[0,\infty)\rightarrow\mathbb{R}$ of class $C^{2}$ such that (compare (1.3)--(1.5))
\begin{equation}
\begin{cases}
\tilde{a}t-\tilde{b}\le g(t)\le\tilde{A}t+\tilde{B}&\text{for } t \ge 0, \\
g'(0)=0, \quad g''\in L^{\infty}([0,\infty)), & \text{and } g''(t)>0, \ t\ge0
\end{cases}
\end{equation}
holds with constants $\tilde{a}, \tilde{b}, \tilde{A}, \tilde{B}>0$. As domain $\Omega$ we take $[|x|>1]$ and observe that (1.6) now
reads as
\begin{equation}
\text{div}\left( \frac{g'(|\nabla u|)}{|\nabla u|} \nabla u \right) = 0 \quad \text{on } [|x|>1].
\end{equation}
In addition to (1.11) we require
\begin{equation}
\int_{0}^{\infty}tg''(t)dt<\infty
\end{equation}
and observe that (1.13) is true, if $g''(t)\le c(1+t)^{-\mu}$ with exponent $\mu>2$. We then introduce
the $g$-catenoids
\begin{equation}
k_{g}^{\pm}(x):=\pm\int_{1}^{|x|}(g')^{-1}\left(\frac{1}{t}\right)dt, \quad |x|\ge1
\end{equation}
(compare formula (2.5) in [8]), which for the choice $g(t):=\sqrt{1+t^{2}}$ take the usual form. Note
that by [8], Lemma 2, the $g$-catenoids are solutions of (1.12) and in addition satisfy ($u=k_{g}^{\pm}$)
\begin{equation}
u\in C^{0}([|x|\ge1]), \quad u(x)=0 \text{ for } |x|=1.
\end{equation}

After these preparations we can deduce from Theorem 1.1:
\begin{corollary}
Let (1.11), (1.13) hold for the function $g$ and define $k_{g}^{\pm}$ as in (1.14). Suppose
that $u\in C^{1}([|x|>1])$ is a weak solution of (1.12) with the property (1.15). Then, if either
$u-k_{g}^{+}$ or $u-k_{g}^{-}$ stays bounded, $u$ is equal to $k_{g}^{+}$ in the first case, whereas
$u=k_{g}^{-}$ holds in the second case.
\end{corollary}

Next we present a result concerning solutions of equation (1.6) whose derivatives in all
coordinate directions vanish on $\partial\Omega$ in a rather strong sense.

\begin{theorem}
Let (1.2), (1.3) and (1.5) hold, whereas condition (1.4) is replaced by
\begin{equation}
|D^{2}f(p)|\le c\frac{1}{1+|p|}, \quad p\in\mathbb{R}^{2N},
\end{equation}
$c$ denoting a positive constant. Let $u\in C^{2}(\Omega,\mathbb{R}^{N})\cap C^{1}(\overline{\Omega},\mathbb{R}^{N})$ denote a solution of (1.6) with
the property
\begin{equation}
|\nabla u(x)|\le c\,\text{dist}(x,\partial\Omega)^{1/2}
\end{equation}
for $x\in\overline{\Omega}$ close to $\partial\Omega$ for some other constant $c>0$ (which for instance follows from
$\partial_{i}u\in C^{0,1/2}$ near $\partial\Omega$ together with $\partial_{i}u=0$ on $\partial\Omega$, $i=1,2$). Then, if in addition $u$ satisfies
\begin{equation}
\lim_{|x|\rightarrow\infty}\frac{|u(x)|}{|x|}<\infty,
\end{equation}
the constancy of $u$ follows.
\end{theorem}

In order to proceed we discuss the situation for the case, when the domain $\Omega$ from (1.2) is
replaced by a half-plane $H$, for instance we may consider
\begin{equation}
H:=[x_{1}>0]:=\{x=(x_{1},x_{2})\in\mathbb{R}^{2}:x_{1}>0\}.
\end{equation}
Concerning problem (1.11) with $H$ in place of $\Omega$ the question of unique solvability has been
analyzed by various authors, and without being complete we mention the papers of Hwang [9],
[10] and of Collin and Krust [11] for further references. It turns out that two solutions $u, v$ must
coincide, if $u-v$ satisfies certain growth restrictions as $|x|\rightarrow\infty$. Concerning the question of
the existence of solutions to (1.1) for the half-plane $H$ the interested should consult the paper
of Jiang, Wang and Zhu [12].

Let us look at equation (1.6) for the half-plane. It is easy to prove:

\begin{theorem}
The results of Theorem 1.1 and of Theorem 1.4 remain true, if the domain $\Omega$ from (1.2) is replaced by the half-plane $H$.
\end{theorem}

For the half-plane $H$ condition (1.17) seems to be somewhat artificial, since the gradient is
required to vanish on $\partial H$. We finish the introduction with a more natural result.

\begin{theorem}
Let $f$ satisfy (1.3), (1.5), (1.16) and consider $H$ as in (1.19). Suppose that
$u\in C^{2}(H,\mathbb{R}^{N})\cap C^{0}(\overline{H},\mathbb{R}^{N})$ is a solution of (1.7) on $H$ such that for a positive constant $c$
\begin{equation}
|\partial_{1}u(x)|\le c\,\text{dist}(x,\partial H) = c\sqrt{x_{1}}
\end{equation}
holds for $x$ near $\partial H=[x_{1}=0]$ (which is true, if $\partial_{1}u=0$ on $\partial H$ and $\partial_{1}u\in C^{0,1/2}$ near $\partial H$).
Then the validity of
\begin{equation}
\lim_{|x|\rightarrow\infty}\frac{|u(x)|}{|x|}<\infty
\end{equation}
implies the existence of points $\alpha, \beta\in\mathbb{R}^{N}$ with the property that $u(x)=\alpha x_{2}+\beta$, $x\in H$.
\end{theorem}

Our paper is organized as follows: Section 2 contains the proofs of Theorem 1.1 and 1.4. In
Section 3 we discuss the half-plane case presenting the proofs of Theorem 1.5 and 1.6.

\section{Proofs of Theorem 1.1 and 1.4}

During the proof of Theorem 1.1 we will benefit from:

\begin{lemma}
Let the density $f$ satisfy (1.3)--(1.5), let $G$ denote a bounded Lipschitz region in
$\mathbb{R}^{2}$ and define $\Omega$ as in (1.2). Consider a function $u$ from the space $C^{1}(\Omega,\mathbb{R}^{N})\cap C^{0}(\overline{\Omega},\mathbb{R}^{N})$
satisfying equation (1.7). Then it holds for any radius $R>0$ such that $\overline{G}\subset[|x|<R]$
\begin{equation}
u\in W^{1,1}([|x|<R]-\overline{G},\mathbb{R}^{N})
\end{equation}
and the Sobolev trace of $u$ on $\partial G=\partial\Omega$ coincides with the classical pointwise trace.
\end{lemma}

\begin{remark}
In the situation of Theorem 1.1 we can apply Lemma 2.1 to the functions $u$ and
$v$ with the result that
\begin{equation}
\Phi(u-v)\in\mathring{W}^{1,1}([|x|<R]-\overline{G},\mathbb{R}^{N})
\end{equation}
holds for arbitrary $\Phi\in C_{0}^{1}([|x|<R])$.
\end{remark}

\begin{proof}[Proof of Lemma 2.1]
We assume for notational simplicity that $N=1$ and consider the case
$G=B_{1}(0)=[|x|<1]$. Fix a radius $R>1$ and choose a small number $\epsilon>0$. We let
$\eta\in C_{0}^{1}([|x|<2R])$ such that $0\le\eta\le1$ and in addition
\begin{equation}
\begin{cases}
\eta=1 & \text{on } [1+\epsilon\le|x|\le R],\\
\eta=0 & \text{on } [|x|\le1+\epsilon/2],\\
|\nabla\eta|\le c/\epsilon & \text{on } [1+\epsilon/2\le|x|\le1+\epsilon],
\end{cases}
\end{equation}
$c$ denoting a positive constant independent of $\epsilon$. Noting that $\eta^{2}u$ belongs to $C_{0}^{1}(\Omega)$ with $\Omega=$
$\mathbb{R}^{2}-\overline{G}=[|x|>1]$ (recall (1.2)), we deduce from (1.7) combined with (2.3)
\begin{equation}
\int_{[1+\epsilon/2<|x|<2R]}\nabla f(\nabla u)\cdot\nabla(\eta^{2}u)dx=0.
\end{equation}
Recalling the boundedness of $\nabla f$ and using (2.3) one more time, we get from equation (2.4)
\begin{align*}
\int_{[1+\epsilon<|x|<R]}\nabla f(\nabla u)\cdot\nabla udx &\le 2\int_{[1+\epsilon/2<|x|<2R]}|\nabla f(\nabla u)|\eta|\nabla\eta||u|dx \\
&\le c(R)+c\int_{[1+\epsilon/2<|x|<1+\epsilon]}|\nabla\eta|dx
\end{align*}
for a constant $c(R)$ and another positive constant $c$ both being independent of $\epsilon$. Of course we
have also used the boundedness of $u$ on the set $[1\le|x|\le2R]$. This yields (recall the bound for
$|\nabla\eta| $ in (2.3))
\begin{equation}
\int_{[1+\epsilon<|x|<R]}\nabla f(\nabla u)\cdot\nabla udx\le K
\end{equation}
for all sufficiently small $\epsilon>0$ and with some number $K\in(0,\infty)$. Using the convexity of $f$ and
condition (1.3) we finally observe that
\[
\nabla f(\nabla u)\cdot\nabla u\ge f(\nabla u)-f(0)\ge a|\nabla u|-b-f(0),
\]
and therefore (2.5) implies $\nabla u\in L^{1}([1<|x|<R],\mathbb{R}^{2})$, which proves (2.1). That the Sobolev
trace on $\partial\Omega=[|x|=1]$ of the function $u$ equals the classical trace then is a consequence of (2.1)
combined with Theorem 3.6.2 and Theorem 3.6.3 (a) in the textbook [13] of Morrey.
\end{proof}

\begin{remark}
Condition (1.4) is not needed during the proof of Lemma 2.1. Moreover, Lemma
2.1 is true for the n-dimensional case, $n\ge2$.
\end{remark}

\begin{proof}[Proof of Theorem 1.1]
Let us assume that all the hypotheses hold, where for notational simplic-
ity we let $N=1$. We consider radii $R$ such that $\overline{G}\subset[|x|<R]$ and observe that according to
Lemma 2.1 it holds $u, v\in W^{1,1}([|x|<2R]-\overline{G})$ and in addition
\begin{equation}
\phi:=\eta^{2}(u-v)\in\mathring{W}^{1,1}([|x|<2R]-\overline{G}),
\end{equation}
where $\eta\in C_{0}^{1}([|x|<2R])$ satisfies
\begin{equation}
\begin{cases}
\eta=1 \text{ on } [|x|\le R], & 0\le\eta\le1, \\
|\nabla\eta|\le\frac{c}{R}. & 
\end{cases}
\end{equation}
For (2.6) we refer to Remark 2.2, in particular (2.2), with $\Phi$ being replaced by $\eta^{2}$. We have
equation (1.7) for $u$ and $v$ with $\varphi$ being replaced by $\phi$ defined in (2.6), which gives after sub-
tracting the results
\[
\int_{[|x|<2R]-\overline{G}}(\nabla f(\nabla u)-\nabla f(\nabla v))\cdot\nabla(\eta^{2}(u-v))dx=0
\]
and in conclusion
\begin{equation}
\begin{aligned}
\int_{[|x|<2R]-\overline{G}}[\nabla f(\nabla u)-\nabla f(\nabla v)]\cdot(\nabla u-\nabla v)\eta^{2}dx = \\
-2\int_{[|x|<2R]-\overline{G}}[\nabla f(\nabla u)-\nabla f(\nabla v)]\cdot\nabla\eta\eta(u-v)dx.
\end{aligned}
\end{equation}

Consider the parameter dependent bilinear form
\[
B_{t,x}(p,q):=D^{2}f(\nabla v(x)+t(\nabla u(x)-\nabla v(x)))(p,q)
\]
for $p,q\in\mathbb{R}^{2}$ and parameters $x, t$, and rewrite (2.8) in the form
\begin{equation}
\begin{aligned}
\int_{[|x|<2R]-\overline{G}}\int_{0}^{1}B_{t,x}(\nabla u-\nabla v,\nabla u-\nabla v)\eta^{2}dtdx = \\
-2\int_{[|x|<2R]-\overline{G}}\int_{0}^{1}B_{t,x}(\eta[\nabla u-\nabla v],(u-v)\nabla\eta)dtdx.
\end{aligned}
\end{equation}
The Cauchy-Schwarz inequality applied to the positive (recall (1.5)) bilinear form $B_{t,x}$ shows
\[
|B_{t,x}(\eta[\nabla u-\nabla v],(u-v)\nabla\eta)|\le\eta~B_{t,x}(\nabla u-\nabla v,\nabla u-\nabla v)^{1/2}|u-v|B_{t,x}(\nabla\eta,\nabla\eta)^{1/2},
\]
hence by Young's inequality we find that
\[
|B_{t,x}(\eta[\nabla u-\nabla v],(u-v)\nabla\eta)|\le\epsilon\eta^{2}B_{t,x}(\nabla u-\nabla v,\nabla u-\nabla v)+c(\epsilon)(u-v)^{2}B_{t,x}(\nabla\eta,\nabla\eta)
\]
valid for all $\epsilon>0$.
Using the latter inequality on the right-hand side of (2.9), we obtain
after appropriate choice of $\epsilon$
\[
\int_{[|x|<2R]-\overline{G}}\int_{0}^{1}\eta^{2}B_{t,x}(\nabla u-\nabla v,\nabla u-\nabla v)dt~dx\le c\int_{[|x|<2R]-\overline{G}}\int_{0}^{1}(u-v)^{2}B_{t,x}(\nabla\eta,\nabla\eta)dt~dx.
\]
Next we recall (1.4), use (2.7) together with assumption (1.9) or (1.10) in order to get for
sufficiently large $R$
\begin{equation}
\int_{[|x|<R]-\overline{G}}\int_{0}^{1}B_{t,x}(\nabla u-\nabla v,\nabla u-\nabla v)dxdt\le C
\end{equation}
for some constant $C\in(0,\infty)$. With (2.10) we go back to (2.9) and recall that the Cauchy-
Schwarz inequality applied to the right-hand side of (2.9) yields (also using (2.7))
\begin{align*}
\int_{[|x|<R]-\overline{G}}\int_{0}^{1}B_{t,x}(\nabla u-\nabla v,\nabla u-\nabla v)dtdx \le \\
2\int_{[|x|<2R]-\overline{G}}\int_{0}^{1}B_{t,x}(\nabla u-\nabla v,\nabla u-\nabla v)^{1/2}|u-v|B_{t,x}(\nabla\eta,\nabla\eta)^{1/2}dtdx,
\end{align*}
hence by Hölder's inequality (and using $\text{spt}(\nabla\eta)\subset[R<|x|<2R]$)
\begin{equation}
\begin{aligned}
&\int_{[|x|<R]-\overline{G}}\int_{0}^{1}B_{t,x}(\nabla u-\nabla v,\nabla u-\nabla v)dtdx \le \\
&\left(\int_{[R<|x|<2R]}\int_{0}^{1}B_{t,x}(\nabla u-\nabla v,\nabla u-\nabla v)dt~dx\right)^{1/2}\cdot\left(\int_{[R<|x|<2R]}\int_{0}^{1}(u-v)^{2}B_{t,x}(\nabla\eta,\nabla\eta)dt~dx\right)^{1/2}.
\end{aligned}
\end{equation}
From (2.10) it follows
\[
\int_{\mathbb{R}^{2}-\overline{G}}\int_{0}^{1}B_{t,x}(\nabla u-\nabla v,\nabla u-\nabla v)dt~dx<\infty,
\]
hence the first term on the right-hand side of (2.11) goes to zero as $R\rightarrow\infty$, whereas the
second integral stays bounded according to the discussion after (2.9), which led us to (2.10).
Altogether we conclude from (2.11) that
\[
\int_{\mathbb{R}^{2}-\overline{G}}\int_{0}^{1}B_{t,x}(\nabla u-\nabla v,\nabla u-\nabla v)dt~dx=0,
\]
hence $\nabla u=\nabla v$ according to (1.6), and (1.8) shows $u=v$, which finishes the proof of The-
orem 1.1.
\end{proof}

For proving Corollary 1.2 we just note that in the case of entire solutions the previous
calculations remain true, if we drop the set $G$ in all formulas.

\begin{proof}[Proof of Theorem 1.4]
For simplicity we let $N=1$ together with $\Omega=[|x|>1]$. It is easy to see
that under the assumptions and with the notation from Theorem 1.4 we can deduce from (2.6)
the validity of (recall the weak form (1.7))
\begin{equation}
\int_{\Omega}D^{2}f(\nabla u)(\nabla\partial_{i}u,\nabla\varphi)dx=0
\end{equation}
for all $\varphi\in C_{0}^{1}(\Omega)$ and $i\in\{1,2\}$ being fixed. Choosing $\eta$ as in the proof of Lemma 2.1 (see (2.3))
and letting $\varphi=\eta^{2}\partial_{i}u$ in equation (2.12), we obtain from (2.12), using the Cauchy-Schwarz
inequality for the bilinear form $D^{2}f(\nabla u)$ followed by an application of Young's inequality, the
estimate
\begin{equation}
\int_{\Omega}\eta^{2}D^{2}f(\nabla u)(\nabla\partial_{i}u,\nabla\partial_{i}u)dx\le c\int_{\Omega}D^{2}f(\nabla u)(\nabla\eta,\nabla\eta)|\partial_{i}u|^{2}dx
\end{equation}
with a uniform constant $c\in(0,\infty)$. Since we assume $u\in C^{1}(\overline{\Omega})$ and according to (2.3), it
follows that we just have to find a bound independent of $\epsilon$ for the quantity
\[
I:=\int_{[1+\epsilon/2<|x|<1+\epsilon]}D^{2}f(\nabla u)(\nabla\eta,\nabla\eta)|\partial_{i}u|^{2}dx.
\]
From assumption (1.16) we get (noting that here (1.4) is sufficient)
\[
I\le\epsilon^{-2}\int_{[1+\epsilon/2\le|x|\le1+\epsilon]}|\nabla u|^{2}dx,
\]
and (1.17) implies $I\le C$ with $C\in(0,\infty)$ independent of $\epsilon$. Referring to (1.5) we find a
number $\lambda(R)\in(0,\infty)$ such that on $[1\le|x|\le R]$ it holds
\[
D^{2}f(\nabla u)(p,p)\ge\lambda(R)|p|^{2}, \quad p\in\mathbb{R}^{2};
\]
the discussion concerning $I$ turns (2.13) into
\begin{equation}
\int_{[1+\epsilon<|x|\le R]}D^{2}f(\nabla u)(\nabla\partial_{i}u,\nabla\partial_{i}u)dx\le c(R)
\end{equation}
for any $\epsilon\ll1$, and (2.14) clearly shows that
\begin{equation}
\partial_{i}u\in W^{1,2}([1<|x|<R])
\end{equation}
for any $R>1$ with vanishing Sobolev trace on $[|x|=1]$ which follows from (1.17) (compare the
remarks at the end of the proof of Lemma 2.1). In order to prove the constancy of $u$, we again
let $R>0$ and choose $\eta\in C_{0}^{1}([|x|<2R])$ with (2.7). The function $\varphi=\eta^{2}\partial_{i}u$ is admissible in
equation (2.12): by the foregoing discussion (recall (2.15) together with $\partial_{i}u=0$ on $[|x|=1]$)
we have that $\eta^{2}\partial_{i}u\in\mathring{W}^{1,2}([1<|x|<2R])$, thus there is a sequence $\varphi_{m}\in C_{0}^{\infty}([1<|x|<2R])$
with the property $\varphi_{m}\rightarrow\varphi$ in the norm of $W^{1,2}([1<|x|<2R])$. Equation (2.12) holds for $\varphi_{m}$
and by passing to the limit $m\rightarrow\infty$ we get the desired result, from which we deduce (compare
(2.8))
\begin{equation}
\int_{[1<|x|<2R]}D^{2}f(\nabla u)(\nabla\partial_{i}u,\nabla\partial_{i}u)\eta^{2}dx=-2\int_{[R<|x|<2R]}D^{2}f(\nabla u)(\eta\nabla\partial_{i}u,\partial_{i}u\nabla\eta)dx
\end{equation}
where we have used the fact that $\text{spt}(\nabla\eta)\subset[R\le|x|\le2R]$. Applying the Cauchy-Schwarz
inequality to $D^{2}f(\nabla u)$ on the right-hand side of (2.16) and then using Hölder's inequality we
get from (2.16)
\begin{equation}
\begin{aligned}
\int_{[1<|x|<2R]}D^{2}f(\nabla u)(\nabla\partial_{i}u,\nabla\partial_{i}u)\eta^{2}dx \le c\left(\int_{[R<|x|<2R]}D^{2}f(\nabla u)(\nabla\partial_{i}u,\nabla\partial_{i}u)\eta^{2}dx\right)^{1/2} \cdot \\
\left(\int_{[R<|x|<2R]}D^{2}f(\nabla u)(\nabla\eta,\nabla\eta)|\partial_{i}u|^{2}dx\right)^{1/2}
\end{aligned}
\end{equation}
with a constant $c\in(0,\infty)$ independent of $R$. We claim
\begin{equation}
\sup_{R>1}R^{-2}\int_{[R<|x|<2R]}|\nabla u|dx<\infty.
\end{equation}
Assuming this for the moment, applying Young's inequality to the right-hand side of (2.17) and
using (1.16), we find that (2.18) implies (passing to $R\rightarrow\infty$)
\begin{equation}
\int_{[|x|>1]}D^{2}f(\nabla u)(\nabla\partial_{i}u,\nabla\partial_{i}u)dx<\infty.
\end{equation}
Now, due to (2.18) and (2.19), the right-hand side of (2.17) actually vanishes as $R\rightarrow\infty$ which
means
\[
\int_{[|x|>1]}D^{2}f(\nabla u)(\nabla\partial_{i}u,\nabla\partial_{i}u)dx=0,
\]
and $\nabla\partial_{i}u=0$ is a consequence of (1.5). This yields the constancy of $\partial_{i}u$, and from $\partial_{i}u=0$ on
$[|x|=1]$ we get $\partial_{i}u=0$ on $\Omega=[|x|>1]$. Since $i$ was arbitrary, $u$ must be a constant function.

It therefore remains to prove (2.18): let $R>2$ and choose $\Psi:\mathbb{R}^{2}\rightarrow[0,1]$ such that $\Psi\in C^{1}$
and $\Psi=1$ on $[R<|x|<2R]$, $\Psi=0$ on $[|x|\le\frac{3}{4}R]\cup[|x|\ge\frac{5}{2}R]$, $|\nabla\Psi|\le c/R$. Then $\Psi^{2}u$ is
admissible in equation (1.7) and with the help of the boundedness of $\nabla f$ we get
\[
\int_{[R<|x|<2R]}\nabla f(\nabla u)\cdot\nabla u~dx\le c\int_{\text{spt}(\nabla\Psi)}|\nabla\Psi||u|dx\le c\int_{[\frac{3}{4}R\le|x|\le\frac{5}{2}R]}R^{-1}|u|dx\le cR^{2},
\]
where the last inequality is a consequence of (1.18). Recalling that by the convexity of $f$
and assumption (1.3) it holds
\[
\nabla f(\nabla u)\cdot\nabla u\ge f(\nabla u)-f(0)\ge a|\nabla u|-b-f(0),
\]
we have shown (2.18), which completes the proof of Theorem 1.4.
\end{proof}

\begin{remark}
If $G$ is a bounded domain in $\mathbb{R}^{3}$, then an inspection of the proof shows, that the
result of Theorem 1.4 still holds, if we replace (1.18) by the stronger condition that $u$ is a bounded
function.
\end{remark}

\section{Proofs of Theorem 1.5 and 1.6}

Just for notational simplicity we consider the scalar case $N=1$ and recall definition (1.19) of the
half-plane $H$ remarking that the particular choice $H=[x_{1}>0]$ just enters through condition
(1.20) from Theorem 1.6 and can easily be adjusted to other choices for $H$.

\begin{itemize}
\item[i)] Letting $\Omega=H$ in Theorem 1.1, the extension of Theorem 1.1 to the half-plane is based
on an appropriate version of Lemma 2.1. We now claim
\begin{equation}
u\in W^{1,1}(H\cap B_{R}(0)), \quad B_{R}(0):=[|x|<R],
\end{equation}
for any $R\in(0,\infty)$. For proving (3.1) we define ($\epsilon>0$ sufficiently small)
\[
\eta(x):=\eta_{1}(x_{1})\eta_{2}(x_{2}), \quad x\in\mathbb{R}^{2},
\]
with functions $\eta_{i}\in C_{0}^{1}(\mathbb{R})$, $0\le\eta_{i}\le1$, $i=1,2$, such that
\begin{equation}
\begin{cases}
\text{spt}(\eta_{2})\subset[-2R, 2R], \quad \eta_{2}=1 \text{ on } [-R,R], \quad |\eta_{2}'|\le c/R, \\
\text{spt}(\eta_{1})\subset[\frac{\epsilon}{2},2R], \quad \eta_{1}=1 \text{ on } [\epsilon,R], \quad |\eta_{1}'|\le c/R \text{ on } [R, 2R], \quad |\eta_{1}'|\le c/\epsilon \text{ on } [\frac{\epsilon}{2}, \epsilon].
\end{cases}
\end{equation}
From (3.2) we get $\eta(x)=1$ if $x\in[\epsilon,R]\times[-R,R]$, and $\eta(x)\ne0$ implies $x\in[\frac{\epsilon}{2},2R]\times[-2R,2R]$,
hence $\text{spt}(\eta)$ is a compact subset of the open half-plane $H$. Equation (2.4) is replaced by
\[
\int_{[\frac{\epsilon}{2},2R]\times[-2R,2R]}\nabla f(\nabla u)\cdot\nabla(\eta^{2}u)dx=0
\]
and following the lines after (2.4), we see that it remains to discuss the term
\[
\int_{[\frac{\epsilon}{2},\epsilon]\times[-2R,2R]}|\eta_{1}'(x_{1})|dx
\]
which however due to the bound imposed on $\eta_{1}'$ (see (3.2)) stays bounded independent of $\epsilon$.
Thus we get $u\in W^{1,1}((0,R)\times(-R,R))$ and in consequence (3.1) is true.

\item[ii)] From (3.1) it follows as before that the Sobolev trace of $u-v$ on $\partial H$ vanishes, if we are
in the situation of Theorem 1.1 with $\Omega=H$ and we deduce
\begin{equation}
\phi:=\eta^{2}(u-v)\in\mathring{W}^{1,1}(H\cap B_{2R}(0)),
\end{equation}
where $R>0$ and $\eta$ now is defined as in (2.7). We have equation (1.6) for the functions $u$ and
$v$, which on account of (3.3) gives
\begin{equation}
\int_{H\cap[|x|<2R]}(\nabla f(\nabla u)-\nabla f(\nabla v))\cdot\nabla\phi dx=0.
\end{equation}
From (3.4) we obtain the half-plane version of (2.8), where one just replaces the domain of
integration by the set $H\cap[|x|<2R]$ and with obvious changes the proof of Theorem 1.1 can
be finished as before.

\item[iii)] If we let $\Omega=H$ in the situation of Theorem 1.4, we first note that (2.15) has to be re-
placed by ($i=1,2$)
\begin{equation}
\partial_{i}u\in W^{1,2}(H\cap[|x|<R]).
\end{equation}
For proving (3.5) we start from equation (2.12) by letting $\varphi=\eta^{2}\partial_{i}u$, but now with the choice
$\eta(x):=\eta_{1}(x_{1})\eta_{2}(x_{2})$ with $\eta_{1}, \eta_{2}$ as in (3.2). The subsequent calculations are easily adjusted.
This also concerns the second part of the proof of Theorem 1.4: we have equation (2.16) with
$H\cap[|x|<2R]$ as domain of integration and with $\eta$ exactly as in (2.7). Proceeding as in Theorem
1.4, we obtain the constancy of $u$.

\item[iv)] We finally come to the proof of Theorem 1.6. We first observe that on account of (1.20) we
obtain (3.5) just for the index $i=1$. For details we refer to the derivation of (2.14) (with $i=1$)
noting that quantity $I$ introduced after (2.13) involves $|\partial_{1}u|^{2}$, to which we can apply (1.20). Of
course, in (2.13) $\eta$ is defined according to (3.2). We then proceed as done after (2.15), letting
$i=1$ and applying the adjustments from iii), and get that $\nabla\partial_{1}u=0$, hence $\partial_{1}u=\alpha$ for some
$\alpha\in\mathbb{R}$. However, (1.20) shows $\alpha=0$ and $u$ just depends on $x_{2}$. From equation (1.6) it follows
directly that we finally arrive at the claim of Theorem 1.6, for details we refer to the discussion
after formula (2.6) of [14].
\end{itemize}
\qed

\section*{References}
\begin{enumerate}
    \item R. Osserman, \textit{A survey of minimal surfaces}. Second Edition, Dover Publications, New York 1986.
    \item E. Kuwert, \textit{A bound for minimal graphs with a normal at infinity}. Calc. Var. Partial Differential Equations 1 no.4 (1993), 407-416.
    \item N. Kutev, F. Tomi, \textit{Existence and nonexistence for the exterior Dirichlet problem for the minimal surface equation in the plane}. Differential Integral Equations 11 no.6 (1998), 917-928.
    \item N. do Espírito-Santo, J. Ripoll, \textit{Some existence results on the exterior Dirichlet problem for the minimal hypersurface equation}. Ann. Inst. H. Poincaré C Anal. Non Linéaire 28 no. 3 (2011), 385-393.
    \item J. Ripoll, F. Tomi, \textit{On solutions to the exterior Dirichlet problem for the minimal surface equation with catenoidal ends}. Adv. Calc. Var. 7 no. 2 (2014), 205-226.
    \item A. Aiolfi, C. Assmann, J. Ripoll, \textit{The Dirichlet problem for the minimal surface equation on unbounded helicoidal domains of $\mathbb{R}^{3}$}. arXiv 2023.
    \item B. Dacorogna, \textit{Direct methods in the calculus of variations}. Applied Mathematical Sciences 78, Springer Verlag, Berlin 1989.
    \item M. Bildhauer, M. Fuchs, \textit{An extension of a theorem of Bers and Finn on the removability of isolated singularities to the Euler-Lagrange equations related to general linear growth problems}. Calc. Var. Partial Differential Equations 61 no. 3 (2022), No. 99, 8 pp.
    \item J.-F. Hwang, \textit{Comparison principles and Liouville theorems for prescribed mean curvature equations in unbounded domains}. Ann. Scuola Norm. Sup. Pisa Cl. Sci 15 no.4 (1988), 341-355.
    \item J.-F. Hwang, \textit{A uniqueness theorem for the minimal surface equation}. Pacific J. Math. 176 no.2 (1996), 357-364.
    \item P. Collin, R. Krust, \textit{Le problème de Dirichlet pour l'équation de surfaces minimales sur des domaines non bornés}. Bull. Soc. Math. France 119 no.4 (1991), 443-462.
    \item G. Jiang, Z. Wang, J. Zhu, \textit{Stability of Edelen-Wang's Bernstein type theorem for the minimal surface equation}. J. of Functional Analysis 284 (2023), 109821.
    \item C.B. Morrey, \textit{Multiple integrals in the calculus of variations}. Springer Verlag, Berlin 2008.
    \item M. Bildhauer, M. Fuchs, \textit{Variants of Bernstein's theorem for variational integrals with linear and nearly linear growth}. Ric. Mat. 73 no.5 (2024), 2911-2923.
\end{enumerate}

\newpage
\noindent
\textbf{Michael Bildhauer} \\
Saarland University \\
Department of Mathematics \\
P.O. Box 15 11 50 \\
66041 Saarbrücken, Germany \\
\texttt{bibi@math.uni-sb.de} \\[0.5cm]
\textbf{Martin Fuchs} \\
Saarland University \\
Department of Mathematics \\
P.O. Box 15 11 50 \\
66041 Saarbrücken, Germany \\
\texttt{fuchs@math.uni-sb.de}

\end{document}